\documentclass[11pt]{article}
\usepackage{amsfonts,amsmath,amsthm,latexsym,color,epsfig}
\usepackage{here}
\usepackage{mathabx}
\usepackage{picins}

\theoremstyle{plain}
\newtheorem{thm}{Theorem}[section]
\newtheorem{fact}[thm]{Fact}
\newtheorem{prop}[thm]{Proposition}
\newtheorem{clm}[thm]{Claim}
\newtheorem{cor}[thm]{Corollary}

\newtheorem{lemma}[thm]{Lemma}

\theoremstyle{definition}
\newtheorem{rem}[thm]{Remark}
\newtheorem{dfn}[thm]{Definition}
\newtheorem{exmp}[thm]{Example}
\newtheorem{quest}[thm]{Question}

\input{epsf}

\def\RR{\mathbb{R}}
\def\qed{\ifvmode\mbox{ }\else\unskip\fi\hskip 1em plus 10fill$\Box$}
\def\eps{\varepsilon}
\def\diam{\mathrm{diam}}

\long\def\ignore#1{}

\def\rc{\advance\leftskip by 0pt plus 40em\rightskip=\leftskip
  \parfillskip=0pt \spaceskip=.3333em \xspaceskip=.5em
  \pretolerance=9999 \tolerance=9999 \hyphenpenalty=9999
  \exhyphenpenalty=9999}

\normalsize
\begin{document}

\title{Borsuk and Ramsey type questions in Euclidean space}
\author{Peter Frankl
\thanks{R\'enyi Institute, Hun\-gar\-ian Academy of Sci\-ences, H--1364 Budapest, POB 127, Hungary. Email: {\tt peter.frankl@gmail.com}. A part of this work was carried out while the author was visiting EPFL in May 2015.}
\and
J\'anos Pach
\thanks{R\'enyi Institute and EPFL, Station 8, CH--1014 Lausanne, Switzerland.
    Email: {\tt pach@cims.nyu.edu}. Supported by Swiss National Science Foundation Grants 200020-144531 and 200021-137574.}
\and
Christian Reiher
\thanks{Fachbereich Math\-e\-matik, Uni\-ve\-rsit\"at Ham\-burg,
Bundes\-stra\ss{}e~55, D-20146 Ham\-burg, Germany,
Email: {\tt Christian.Reiher@uni-hamburg.de}}
\and
Vojt\v{e}ch  R\"odl\thanks{Department of Mathematics, Emory University, Atlanta, GA 30322, USA, Email: {\tt rodl@mathcs.emory.edu}. Supported by NSF grants DMS-1301698 and DMS-1102086.}}

\date{}

\maketitle

\centerline{\em Dedicated to Ron Graham on the occasion of his $80$th birthday}

\begin{abstract}
We give a short survey of problems and results on (1) diameter graphs and hypergraphs, and (2) geometric Ramsey theory. We also make some modest contributions to both areas. Extending a well known theorem of Kahn and Kalai which disproved Borsuk's conjecture, we show that for any integer $r\ge 2$, there exist $\varepsilon=\varepsilon(r)>0$ and $d_0=d_0(r)$ with the following property. For every $d\ge d_0$, there is a finite point set $P\subset\mathbb{R}^d$ of diameter $1$ such that no matter how we color the elements of $P$ with fewer than $(1+\varepsilon)^{\sqrt{d}}$ colors, we can always find $r$ points of the same color, any two of which are at distance $1$.
\smallskip

Erd\H{o}s, Graham, Montgomery, Rothschild, Spencer, and Strauss called a finite point set ${P\subset\mathbb{R}^d}$ {\em Ramsey} if for every $r\ge 2$, there exists a set $R=R(P,r)\subset\mathbb{R}^D$ for some $D\ge d$ such that no matter how we color all of its points with $r$ colors, we can always find a monochromatic congruent copy of $P$. If such a set $R$ exists with the additional property that its diameter is the same as the diameter of $P$, then we call $P$ {\em diameter-Ramsey}. We prove that, in contrast to the original Ramsey property, (a) the condition that $P$ is diameter-Ramsey is not hereditary, and (b) not all triangles are diameter-Ramsey. We raise several open questions related to this new concept.

\end{abstract}

\section{Introduction}

The aim of this article is twofold. In the spirit of Graham-Yao \cite{GrY90}, we give a ``whirlwind tour'' of two areas of Geometric Ramsey Theory, and make some modest contributions to them.
\smallskip

The {\em diameter} of a finite point set $P$, denoted by ${\rm diam}(P)$, is the largest distance that occurs between two points of $P$. Borsuk's famous conjecture~\cite{Bor33}, restricted to finite point sets, states that any such set of unit diameter in $\mathbb{R}^d$ can be colored by $d+1$ colors so that no two points of the same color are at distance {\em one}. This conjecture was disproved in a celebrated paper of Kahn and Kalai~\cite{KaK93}. We extend the theorem of Kahn and Kalai as follows.

\medskip
\noindent{\bf Theorem 1.} {\em For any integer $r\ge 2$, there exist $\varepsilon=\varepsilon(r)>0$ and $d_0=d_0(r)$ with the following property. For every $d\ge d_0$, there is a finite point set $P\subset\mathbb{R}^d$ of diameter $1$ such that no matter how we color the elements of $P$ with fewer than $(1+\varepsilon)^{\sqrt{d}}$ colors, we can always find $r$ points of the same color, any two of which are at distance $1$.}
\medskip

In a seminal paper of Erd\H{o}s, Graham, Montgomery, Rothschild, Spencer, and Strauss~\cite{ErGM73}, the following notion was introduced.  A finite set $P$ of points in a Euclidean space is a {\em Ramsey configuration} or, briefly, is {\em Ramsey}
if for every $r\ge 2$, there exists an integer $d=d(P,r)$ such that no matter how we color all points of~$\mathbb{R}^d$ with $r$ colors, we can always find a monochromatic subset of $\mathbb{R}^d$ that is congruent to $P$. In two follow-up articles~\cite{ErGM75a},~\cite{ErGM75b}, Erd\H os, Graham, and their coauthors established many important properties of these sets.

In the present paper, we introduce a related notion.

\medskip
\noindent{\bf Definition 2.} {\em A finite set $P$ of points in a Euclidean space is {\em diameter-Ramsey} if for every integer $r\ge 2$, there exist an integer $d=d(P,r)$ and a finite subset $R\subset \mathbb{R}^d$ with ${\rm diam}(R)={\rm diam}(P)$ such that no matter how we color all points of $R$ with $r$ colors, we can always find a monochromatic subset of $R$ that is congruent to $P$.}
\medskip

Obviously, every diameter-Ramsey set is Ramsey, but the converse is not true. For example, we know that all triangles are Ramsey, but not all of them are diameter-Ramsey.

\medskip
\noindent{\bf Theorem 3.} {\em All acute and all right-angled triangles are diameter-Ramsey.}

\medskip
\noindent{\bf Theorem 4.} {\em No triangle that has an angle larger than $150^\circ$ is diameter-Ramsey.}
\medskip

There is another big difference between the two notions: By definition, every subset of a Ramsey configuration is Ramsey. This is not the case for diameter-Ramsey sets.

\medskip
\noindent{\bf Theorem 5.} {\em The $7$-element set consisting of a vertex of a $6$-dimensional cube and its $6$ adjacent vertices is not diameter-Ramsey.}
\medskip

We will see that the vertex set of a cube (in fact, the vertex set of any brick) is diameter-Ramsey; see Lemma~\ref{brick}. Therefore, the property that a set is diameter-Ramsey is not hereditary.
\smallskip

It appears to be a formidable task to characterize all diameter-Ramsey simplices. It easily follows from the definition that all regular simplices are diameter-Ramsey; see Proposition~\ref{simplex}. We will show that the same is true for ``almost regular'' simplices.
\medskip

\noindent{\bf Theorem 6.}
{\em For every integer $n\ge 2$, there exists a positive real number $\eps=\eps(n)$ such that every $n$-vertex simplex whose side lengths belong to the interval
$[1-\eps, 1+\eps]$ is diameter-Ramsey.}
\medskip

This article is organized as follows: In Section 2, we give a short survey of problems and results on the structure of diameters and related coloring questions. In Section 3, we prove Theorem 1. In Section 4, we establish some simple properties of diameter-Ramsey sets and prove Theorems 3, 4, and 6, in a slightly stronger form. The proof of Theorem 5 is presented in Section 5. The last section contains a few open problems and concluding remarks.

\section{A short history}

{\bf I. The number of edges of diameter graphs and hypergraphs.}
Hopf and Pannwitz~\cite{HoP34} noticed that in any set $P$ of $n$ points in
the plane, the diameter occurs at most~$n$ times. In other words, among the
${n\choose 2}$ distances between pairs of points from~$P$ at most $n$ are
equal to $\diam(P)$. This bound can be attained for every $n\ge 3$. For odd $n$
this is shown by the vertex set of a regular $n$-gon, and for even $n$ it is
not hard to observe that one may add a further point to the vertex set of a
regular $(n-1)$-gon so as to obtain such an example. In fact all extremal
configurations were characterized by Woodall~\cite{Wo71}.

\smallskip

The same question in $\mathbb{R}^3$ was raised by V\'azsonyi, who conjectured
that the maximum number of times the diameter can occur among $n\ge 4$ points
in $3$-space is $2n-2$. V\'azsonyi's conjecture was proved independently
by Gr\"unbaum~\cite{Gr56}, by Heppes~\cite{He56}, and by Straszewicz~\cite{St57};
see also~\cite{Sw08} for a simple proof.
The extremal configurations were characterized in terms of ball polytopes
by Kupitz, Martini, and Perles~\cite{KuMP10}.

\smallskip

In dimensions larger than $3$, the nature of the problem is radically different.

\begin{thm}\label{erdos}\rm (Erd\H os~\cite{Er60}) For any integer $d>3$, the maximum number of occurrences
of the diameter (and, in fact, of any fixed distance) in a set of $n$ points in $\mathbb{R}^d$
is $\frac12\left(1-\frac{1}{\lfloor d/2\rfloor}+o(1)\right)n^2.$
\end{thm}

More recently, Swanepoel~\cite{Sw09} determined the exact maximum number of
appearances of the diameters for all $d>3$ and all $n$ that are sufficiently
large depending on $d$.

\smallskip

The {\em diameter graph} associated with a set of points $P$ is a graph with
vertex set $P$, in which two points are connected by an edge if and only if
their distance is $\diam(P)$. Erd\H{o}s noticed that there is an intimate
relationship between the above estimates for the number of edges of diameter
graphs and the following attractive conjecture of Borsuk~\cite{Bor33}:
Every (finite) $d$-dimensional point set can be decomposed into at most $d+1$
sets of smaller diameter. If it were true, this bound would be best possible, as
demonstrated by the vertex set of a regular simplex in $\mathbb{R}^d$.

\smallskip

One can generalize the notion of diameter graph as follows. Given a point set
$P\subset \mathbb{R}^d$ and an integer $r\ge 2$, let $H_r(P)$ denote the hypergraph
with vertex set $P$ whose hyperedges are all $r$-element subsets
$\{p_1,\ldots,p_r\}\subseteq P$ with $|p_i-p_j|=\diam(P)$
whenever $1\le i\neq j\le r$. Obviously, $H_2(P)$ is the diameter graph of $P$,
and $H_r(P)$ consists of the vertex sets of all {\em $r$-cliques}
(complete subgraphs with $r$ vertices) in the diameter graph.
Note that every $r$-clique corresponds to a regular $(r-1)$-dimensional simplex
with side length ${\rm diam}(P)$. We call $H_r(P)$ the {\em $r$-uniform diameter hypergraph} of $P$.

\smallskip

It was conjectured by Schur that the Hopf-Pannwitz theorem mentioned at the
beginning of this subsection can be extended to higher dimensions in the following
way: For any $d\ge 2$ and any $d$-dimensional $n$-element point set $P$,
the hypergraph $H_d(P)$ has at most $n$ hyperedges.
This was proved for $d=3$ by Schur, Perles, Martini, and Kupitz~\cite{ScPMK03}.
Building on work of Mori\'c and Pach~\cite{MoP15}, the case $d=4$ was resolved
by Kupavskii~\cite{Ku14}, and the general case of Schur's conjecture was subsequently settled by Kupavskii and Polyanskii~\cite{KuP14}.

\smallskip

However, for $2<r<d$ we know very little about the number of edges of the
diameter hypergraphs $H_r(P)$ and it would be interesting to investigate
this matter further.
\medskip

\noindent{\bf II. The chromatic number of diameter graphs and hypergraphs.} Erd\H os~\cite{Er46} pointed out that if we could prove that the number of edges
of the diameter graph of every $n$-element point set $P\subset\mathbb{R}^d$ is
smaller than $\frac{d+1}{2}n$, then this would imply that there is a vertex of
degree at most $d$. Hence, the chromatic number of the diameter graph would be at
most $d+1$, and the color classes of any proper coloring with $d+1$ colors would
define a decomposition of $P$ into at most $d+1$ pieces of smaller diameter, as
required by Borsuk's conjecture. For $d=2$ and $3$, this is the case. However, as is shown by Theorem~\ref{erdos}, in higher dimensions the number of edges of an $n$-vertex diameter graph can grow quadratically in $n$. Based on this, Erd\H os later suspected that Borsuk's conjecture may be false (personal communication).
This was verified only in 1993 by Kahn and Kalai~\cite{KaK93}.

\smallskip

Using a theorem of Frankl and Wilson~\cite{FrW81}, Kahn and Kalai established
the following much stronger statement.

\begin{thm}
{\rm (Kahn-Kalai)} For any sufficiently large $d$,
there is a finite point set $P$ in the \mbox{$d$-dimensional} Euclidean space such that
no matter how we partition it into fewer than $(1.2)^{\sqrt{d}}$ parts, at least
one of the parts contains two points whose distance is $\diam(P)$.
\end{thm}

In other words, the chromatic number of the diameter graph of $P$ is at least
$(1.2)^{\sqrt{d}}$.
Today Borsuk's conjecture is known to be false for all dimensions
$d\ge 64$; cf.~\cite{JeB14}.

\smallskip

\begin{dfn} The {\em chromatic number} of a hypergraph $H$ is the smallest number
$\chi=\chi(H)$ with the property that the vertex set of $H$ can be colored
with $\chi$ colors such that no hyperedge of $H$ is monochromatic.
\end{dfn}

Clearly, we have
\[
\chi(H_r(P))\le \chi(H_{r-1}(P))\le\ldots \le\chi(H_2(P))\,,
\]
for every $P$ and $r\ge 2$. Moreover,
\[
\chi(H_r(P))\le \biggl\lceil\frac{\chi(H_2(P))}{r-1}\biggr\rceil\,.
\]
To see this, take a proper coloring of the diameter graph $H_2(P)$ with the
minimum number, ${\chi=\chi(H_2(P))}$, of colors and let $P_1,\ldots, P_{\chi}$
be the corresponding color classes. Coloring all elements of
\[
P_{(i-1)(r-1)+1}\cup P_{(i-1)(r-1)+2}\cup\ldots\cup P_{i(r-1)}
\]
with color $i$ for $1\le i\le\frac{\chi}{r-1}$, we obtain a proper coloring of the hypergraph $H_r(P)$. (Here we set $P_s=\emptyset$ for all $s>\chi$.)

\smallskip

Using the above notation, the Kahn-Kalai theorem states that for any
sufficiently large integer~$d$, there exists a set $P\subset\mathbb{R}^d$
with $\chi(H_2(P))\ge (1.2)^{\sqrt{d}}$. According to a result of
Schramm~\cite{Sch88}, we have
$\chi(H_2(P))\le \bigl(\sqrt{3/2}+\varepsilon\bigr)^d$
for every $\varepsilon>0$, provided that $d$ is sufficiently large.

\smallskip

In the next section, we prove Theorem 1 stated in the Introduction. It extends the Kahn-Kalai theorem to
$r$-uniform diameter hypergraphs with $r\ge 2$. Using the above notation, we will prove the following.

\begin{thm}\label{thm:KK}
For any integer $r\ge 2$, there exist $\varepsilon=\varepsilon(r)>0$ and
$d_0=d_0(r)$ with the following property. For every $d\ge d_0$, there is a
finite point set $P\subset\mathbb{R}^d$ of diameter $1$ such that
\[
\chi(H_r(P))\ge (1+\varepsilon)^{\sqrt{d}}\,.
\]
That is, for any partition of $P$ into fewer than $(1+\varepsilon)^{\sqrt{d}}$ parts at least one of the parts contains $r$ points any two of which are at distance $1$.
\end{thm}

\medskip

\noindent{\bf III. Geometric Ramsey theory.} Recall from the Introduction that, according to the definition of Erd\H os, Graham {\it et al.}~\cite{ErGM73}, a finite set of points in some Euclidean space is said to be {\em Ramsey} if for every $r\ge 2$, there exists an integer $d=d(P,r)$ such that no matter how we color all points of $\mathbb{R}^d$ with $r$ colors, we can always find a monochromatic subset of $\mathbb{R}^d$ that is congruent to~$P$. Erd\H os, Graham {\it et al.} proved, among many other results, that every Ramsey set is {\em spherical}, i.e., embeddable into the surface of a sphere. Later Graham~\cite{Gr94} conjectured that the converse is also true: every spherical configuration is {\em Ramsey}. An important special case of this conjecture was settled by Frankl and R\"odl.

\begin{thm}\label{franklrodl} {\rm ~\cite{FrR90}} Every simplex is Ramsey.
\end{thm}

It was shown in \cite{ErGM73} that the class of all Ramsey sets is closed both under taking subsets and taking Cartesian products. This implies

\begin{cor} {\rm \cite{ErGM73}}\label{bricks}
All {\em bricks}, i.e., Cartesian products of finitely many $2$-element sets, are Ramsey.
\end{cor}

Further progress in this area has been rather slow. The first example of a planar Ramsey configuration with at least {\em five} elements was exhibited by K\v{r}\'i\v{z}, who showed that every regular polygon is Ramsey. He also proved that the same is true for every Platonic solid. Actually, he deduced both of these statements from the following more general theorem.

\begin{thm} {\rm ~\cite{Kr91}}
If there is a soluble group of isometries acting on a finite set of points $P$ in $\mathbb{R}^d$, which has at most $2$ orbits, then $P$ is Ramsey.
\end{thm}

Graham's conjecture is still widely open. In fact, it is not even known whether all quadrilaterals inscribed in a circle are Ramsey.
\smallskip

An alternative conjecture has been put forward by Leader, Russell, and Walters~\cite{LRW12}.
They call a point set {\it transitive} if its symmetry group is transitive. A subset of a transitive set is said to be {\it subtransitive}. Leader {\it et al.} conjecture that a set is Ramsey if and only if it is subtransitive. It is not obvious {\em a priori} that this conjecture is different from Graham's, that is, if there exists any spherical set which is not subtransitive. However, this was shown to be the case in~\cite{LRW12}. In~\cite{LRW11}  the same authors showed further that not all quadrilaterals inscribed in a circle are subtransitive.

\smallskip

The ``compactness'' property of the chromatic number, established by Erd\H{o}s and de Bruijn~\cite{BrE51}, implies that for every Ramsey set $P$ and every positive integer $r$, there exists a {\em finite} configuration $R=R(P,r)$ with the property that {\em no matter how we color the points of $R$ with $r$ colors, we can find a congruent copy of $P$ which is monochromatic}. Following the (now standard) notation introduced by Erd\H os and Rado, we abbreviate this property by writing
$$R\longrightarrow (P)_r.$$
In Section 4, we address the problem how small the diameter of such a set $R$ can be. In particular, we investigate the question whether there exists a set $R$ with ${\rm diam}(R)={\rm diam}(P)$ such that $R\longrightarrow (P)_r.$ If such a set exists for every $r$, then according to Definition 2 (in the Introduction), $P$ is called {\rm diameter-Ramsey.}

\section{Proof of Theorem 1}
\label{sec:KK}

The proof of Theorem 1, reformulated as Theorem~\ref{thm:KK}, is based on the construction used by Kahn and Kalai in~\cite{KaK93}.

\smallskip

Suppose for simplicity that $d={2n\choose 2}$ holds for some {\em even}
integer $n$ and set $[2n]=\{1,2,\ldots,2n\}$. The construction takes place in
$\RR^d$ and in the following we will index the coordinates of this space by the
$2$-element subsets of $[2n]$.

\smallskip

To each partition $[2n]=X\cup Y$ of $[2n]$ into two $n$-element subsets
$X$ and $Y$, we assign the point $p(X,Y)=p(Y,X)\in\mathbb{R}^d$ whose
coordinate $p_T(X, Y)$ corresponding to some unordered pair $T\subseteq [2n]$
is given by
\[
p_T(X,Y)=\left\{
  \begin{array}{ll}
  1 & {\mbox {\rm if}}\;\; |T\cap X|=|T\cap Y|=1, \\
  0 & {\mbox {\rm otherwise.}}
  \end{array}
  \right.
\]
Let $P\subseteq\mathbb{R}^d$ be the set of all such points $p(X,Y)$.
We have $|P|=\frac12{2n\choose n}$.
\smallskip

Each point $p(X,Y)\in P$ has precisely $|X|\,|Y|=n^2$ nonzero coordinates. The squared Euclidean distance between $p(X,Y)$ and $p(X',Y')$, for two different partitions of $[2n]$, is equal to the number of coordinates in which $p(X,Y)$ and $p(X',Y')$ differ. The number of coordinates in which both $p(X,Y)$ and $p(X',Y')$ have a $1$ is equal to
\[
|X\cap X'||Y\cap Y'|+|X\cap Y'||X'\cap Y|\,.
\]
Denoting $|X\cap X'|=|Y\cap Y'|$ by $t$, the last expression is equal to $t^2+(n-t)^2$. Thus, we have
\[
\|p(X,Y)-p(X',Y')\|^2=2n^2-2(t^2+(n-t)^2)\,,
\]
which attains its maximum for $t=\frac{n}{2}$. The maximum is $n^2$, so that
$\diam(P)=n.$

\begin{fact}\label{fact:n2}
An $r$-element subset $\{p(X_1,Y_1),\ldots,p(X_r,Y_r)\}\subseteq P$ is a hyperedge of $H_r(P)$, the $r$-uniform diameter hypergraph of $P$, if and only if
\[
|X_i\cap X_j|=\frac{n}{2}\;\;\; \mbox{for all }\, 1\le i\neq j\le r. \;\;\;\;\;\;\; \Box
\]
\end{fact}

We need the following important special case of a result of Frankl and
R\"odl~\cite{FrR87} from extremal set theory.
The set of all $n$-element subsets of $[2n]$ is denoted by ${[2n]\choose n}$.

\begin{thm} {\rm \cite{FrR87}}\label{thm:FR}
For every integer $r\ge 2$, there exists $\gamma=\gamma(r)>0$ with the following
property. Every family of subsets ${\mathcal F}\subseteq {[2n]\choose n}$ with
$|{\mathcal F}|\ge (2-\gamma)^{2n}$ has $r$ members,
$F_1,\ldots,F_r\in{\mathcal F}$, such that
\[
|F_i\cap F_j|=\left\lfloor\frac{n}{2}\right\rfloor\;\;\; \mbox{for all }\, 1\le i\neq j\le r\,.
\]
\end{thm}

\smallskip

To establish Theorem~\ref{thm:KK}, fix a subset $Q$ of the set $P$ defined above.
The elements of $Q$ are points $p(X,Y)\in\mathbb{R}^d$ for certain partitions
$[2n]=X\cup Y$. Let ${\mathcal F}(Q)\subseteq{[2n]\choose n}$ denote the family
of all sets $X$ and $Y$ defining the points in $Q$. Notice that
$|{\mathcal F}(Q)|=2\,|Q|$.

\smallskip

By definition, $\chi=\chi(H_r(P))$ is the smallest number for which there is
a partition
\[
P=Q_1\cup\ldots\cup Q_{\chi}
\]
such that no $Q_k$ contains any hyperedge belonging to $H_r(P)$.
According to Fact~\ref{fact:n2}, this is equivalent to the condition that
${\mathcal F}(Q_k)$ does not contain $r$ members such that any two have precisely
$\frac{n}{2}$ elements in common. Now Theorem~\ref{thm:FR} implies that
\[
|{\mathcal F}(Q_k)|=2\,|Q_k|<\bigl(2-\gamma(r)\bigr)^{2n}
\quad \text{whenever} \quad 1\le k\le\chi\,.
\]
Thus, we have
\[
|P|=\sum_{k=1}^{\chi}|Q_k|<\frac{\chi}{2}\bigl(2-\gamma(r)\bigr)^{2n}\,.
\]
Comparing the last inequality with the equation $|P|=\frac12{2n\choose n}$,
we obtain
\[
\chi=\chi(H_r(P))>\frac{{2n\choose n}}{\bigl(2-\gamma(r)\bigr)^{2n}}
>\left(1+\frac{\gamma(r)}{3}\right)^{\sqrt{2d}}\,.
\]
This completes the proof of Theorem~\ref{thm:KK}.

\medskip

The proof of Theorem~\ref{thm:KK} gives the following result. The {\em regular} simplex $S_r$ with $r$ vertices and {\em unit side length} is not only a Ramsey configuration, but for every $k$ there exists set $P(k)\subseteq \mathbb{R}^d$ of {\em unit diameter} with $d\le c(r)\log^2 k$ such that no matter how we color $P(k)$ with $k$ colors, it contains a monochromatic congruent copy of $S_r$. (Here $c(r)>0$ is a suitable constant that depends only on~$r$.)

\section{Diameter-Ramsey sets -- Proofs of Theorems 3, 4, and 6}

According to Definition 2 (in the Introduction), a finite point set $P$ is diameter-Ramsey if for every $r\ge 2$, there exists a finite set $R$ in some Euclidean space with ${\rm diam}(R)={\rm diam}(P)$ such that no matter how we color all points of $R$ with $r$ colors, we can always find a monochromatic subset of $R$ that is congruent to $P$. Before proving Theorems 3, 4, and 6, we make some general observations about diameter-Ramsey sets.

\medskip
\begin{prop}\label{simplex}
Every regular simplex is diameter-Ramsey.
\end{prop}

\noindent{\bf Proof.} Let $P$ be (the vertex set of) a $d$-dimensional regular simplex. For a fixed integer $r\ge 2$, let $R$ be an $rd$-dimensional regular simplex of the same side length. By the pigeonhole principle, no matter how we color the vertices of $R$ with $r$ colors, at least $d+1$ of them will be of the same color, and they induce a congruent copy of $P$. \hfill $\Box$

\medskip

Recall that a {\em brick} is the vertex set of the Cartesian product of finitely many $2$-element sets.

\begin{lemma}\label{brick}
If $P$ and $Q$ are diameter-Ramsey sets, then so is their Cartesian product $P\times Q$. Consequently, any brick is diameter-Ramsey.
\end{lemma}

\noindent{\bf Proof.}  It was shown in \cite{ErGM73} that for any Ramsey sets $P$ and $Q$, their Cartesian product,
\[
P\times Q = \{ p\times q\;|\;p\in P, q\in Q\}\,,
\]
is also a Ramsey set. Their argument, combined with the equation
\[
\diam^2(P\times Q)=\diam^2(P)+\diam^2(Q)\,,
\]
proves the lemma. \hfill $\Box$
\medskip

\noindent{\bf Proof of Theorem 3.}
Consider a right-angled triangle $T$ whose legs are of length $l_1$ and $l_2$. Let $P$ (resp., $Q$) be a set consisting of two points at distance $l_1$ (resp., $l_2$) from each other, so that we have
$T\subseteq P\times Q$. By Lemma~\ref{brick}, $P\times Q$ is diameter-Ramsey. Since $\diam(T)=\diam(P\times Q)$, we also have that $T$ is diameter-Ramsey.
\smallskip

Now let $T$ be an acute triangle with sides $a$, $b$, and $c$, where $a\le b\le c$.
Set
\[
l_1=\sqrt{c^2-a^2}, \quad l_2=\sqrt{c^2-b^2}, \quad \text{ and } \quad
x=\sqrt{a^2+b^2-c^2}\,.
\]
Since $T$ is acute, we have $a^2+b^2-c^2>0$. Therefore, $x$ is well defined. We have $l_1\ge l_2\ge 0$. Suppose first that $l_1\ge l_2>0$. Let $T_0$ be a right angled triangle with legs $l_1$ and $l_2$, and let $S$ be an equilateral triangle of side length $x$. We have $a^2=l_2^2+x^2$, $b^2=l_1^2+x^2$, and $c^2=l_1^2+l_2^2+x^2$. Thus,
\[
T\subseteq T_0\times S \quad \text{ and } \quad \diam(T)=\diam(T_0\times S)=c\,.
\]
By Proposition~\ref{simplex} and Lemma~\ref{brick}, we conclude that $T$ is diameter-Ramsey. In the remaining case, we have $l_2=0$. Now $T_0$ degenerates into a line segment or a point. It is easy to see that the above proof still applies. \hfill$\Box$
\medskip

We will prove Theorem 4 in a more general form. For this, we need a definition.

\begin{dfn}\label{degenerate}
Let $t$ be a positive integer. A finite set of points $P$ in some Euclidean space is said to be {\em $t$-degenerate} if it has a point $p\in P$ such that for the vertex set $S$ of any regular $t$-dimensional simplex with $p\in S$ and ${\rm diam}(S)={\rm diam}(P)$, we have
      $${\rm diam}(P\cup S)>{\rm diam}(P).$$
\end{dfn}

\begin{thm}\label{degthm}
Let $t\ge 1$ and let $P$ be a finite $t$-degenerate set of points in some Euclidean space, which contains the vertex set of a regular $t$-dimensional simplex of side length ${\rm diam}(P)$. Then $P$ is not diameter-Ramsey.
\end{thm}

\noindent{\bf Proof.} Suppose for contradiction that $P$ is diameter-Ramsey. This implies that there exists a set~$R$ with ${\rm diam}(R)={\rm diam}(P)$ such that no matter how we color it by two colors, it always contains a monochromatic congruent copy of $P$.

Color the points of $R$ with red and blue, as follows. A point is colored {\em red} if it belongs to a subset $S\subset R$ that spans a $t$-dimensional simplex of side length ${\rm diam}(R)
$. Otherwise, we color it {\rm blue}. Let $P'$ be a monochromatic copy of $P$. By the assumptions, $P'$ contains the vertices of a regular $t$-dimensional simplex of side length ${\rm diam}(P)$, and all of these vertices are red. Since $P$ is $t$-degenerate, the point of $P'$ corresponding to $p$ is blue, which is a contradiction. \hfill $\Box$

\medskip

Theorem 4 is an immediate corollary of Theorem~\ref{degthm} and the following statement.

\begin{lemma}\label{150degrees1}
Every triangle that has an angle larger than $150^\circ$ is $1$-degenerate.
\end{lemma}

With no danger of confusion, for any two points $p$ and $p'$, we write $pp'$ to denote both the segment connecting them and its length.
\smallskip

To establish Lemma~\ref{150degrees1}, it is sufficient to verify the following.

\begin{lemma}\label{150degrees2}
Let $T=\{p_1,p_2,p_3\}$ be the vertex set of a triangle and $q$ another point in some Euclidean space such that
$$\max(p_2q, p_3q)\le p_1q\le p_2p_3.$$
Then the angle of $T$ at $p_1$ is at most $150^\circ$.
\end{lemma}

First, we show why Lemma~\ref{150degrees2} implies Lemma~\ref{150degrees1}. Let $T=\{p_1,p_2,p_3\}$ be a triangle whose angle at $p_1$ is larger than $150^\circ$, so that $\diam(T)=p_2p_3$. Suppose without loss of generality that $\diam(T)=1$. To prove that $T$ is $1$-degenerate, it is enough to show that for any unit segment $S=p_1q$, we have $\diam(T\cup S)>1$. Suppose not. Then we have
$\max(p_2q,p_3q)\le p_1q=p_2p_3=1.$ Hence, by Lemma~\ref{150degrees2}, the angle of $T$ at $p_1$ is at most $150^\circ$, which is a contradiction.

\medskip
\noindent{\bf Proof of Lemma~\ref{150degrees2}.} Proceeding indirectly, we assume that
\begin{equation}\label{eq:more150}
\sphericalangle{p_2p_1p_3}>150^\circ\,.
\end{equation}
Let $\Pi$ denote a ($2$-dimensional) plane containing $T$, and let $q'$ denote the orthogonal projection of $q$ to $\Pi$. In the plane $\Pi$, let $g$ and $h$ denote the perpendicular bisectors of the segments $p_1p_2$ and $p_2p_3$, respectively.

\smallskip

\parpic[l]{\includegraphics[width=6cm]{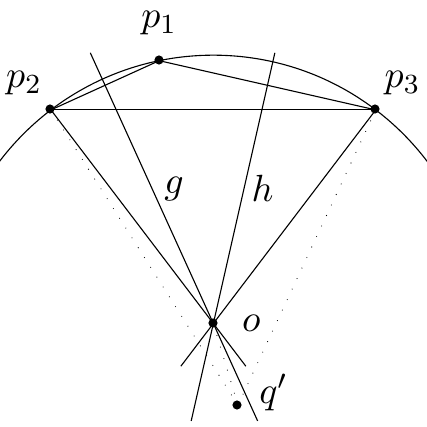}} Since $p_1q\ge p_2q$, we have $p_1q'\ge p_2q'$. Thus, $q'$ belongs to the closed half-plane of $\Pi$ bounded by $g$ where $p_2$ lies. By symmetry, $q'$ belongs to the half-plane bounded by $h$ that contains~$p_3$. This implies that the intersection of these two half-planes is nonempty. In particular, $p_1$ cannot be an interior point of $p_2p_3$ and, by \eqref{eq:more150}, it follows that the triangle $T$ must be non-degenerate. Hence, $g$ and $h$ must meet at a point $o$, the circumcenter of $T$.
\smallskip

Due to the inscribed angle theorem, we have
\[
\sphericalangle{p_2p_1p_3}+\tfrac12 \sphericalangle{p_2op_3}=180^\circ
\]
and hence $\sphericalangle{p_2op_3}< 60^\circ$ by~\eqref{eq:more150}.
This, in turn, implies that $p_2o, p_3o> p_2p_3$.
Thus, we have $$p_2q'\le p_2q\le p_2p_3<p_2o$$ and, in particular, $q'\ne o$. If one side of a triangle is smaller than another, then the same is true for the opposite angles. Applying this to the triangle $p_2q'o$, we obtain that $\sphericalangle{q'op_2}<90^\circ$. Analogously, we have
$\sphericalangle{q'op_3}<90^\circ$, which contradicts the position of $q'$ described in the previous paragraph. \hfill$\Box$

\medskip
We have been unable to answer

\begin{quest}\label{qu:obtuse}
Does there exist any obtuse triangle that is diameter-Ramsey?
\end{quest}

We would like to remark, however, that the answer would be affirmative if we would
just consider colourings with two colours. This is shown by the following example.

\begin{exmp}
Let $R$ be the vertex set of a regular heptagon $p_1p_2\dots p_7$ and let 
$P=\{p_1, p_2, p_4\}$. Clearly, $P$ is the vertex set of an obtuse triangle
having an angle of size $\tfrac 47\cdot 180^\circ>90^\circ$ and ${\diam(R)=\diam(P)}$. 
Moreover, we 
have $R\longrightarrow (P)_2$, because the triple system with vertex set~$R$
whose edges are all sets of the form $\{p_i, p_{i+1}, p_{i+3}\}$ (the addition being 
performed modulo $7$) is known to be isomorphic to the Fano plane, which 
in turn is known to have chromatic number $3$.

\end{exmp}

It seems to be quite difficult to characterize all diameter-Ramsey simplices. According to Proposition~\ref{simplex}, every regular simplex is diameter-Ramsey. Theorem 6 states that this remains true for ``almost regular'' simplices. It is a direct corollary of the following statement.

\begin{lemma}
Every simplex $S$ with vertices $p_1, p_2, \ldots, p_n$ satisfying
\[
\sum_{1\le i<j\le n}(p_ip_j)^2\ge \left(\binom{n}{2}-1\right)\diam^2(S)
\]
is diameter-Ramsey.
\end{lemma}

\noindent{\bf Proof.}
Suppose without loss of generality that $\diam(S)=p_1p_2=1$. Our strategy is to embed~$S$ into the Cartesian product $R$ of $1+\binom{n}{2}$ regular simplices, some of which might degenerate to a point. We will be able to achieve this, while making sure that $\diam(R)=1$. Thus, in view of Proposition~\ref{simplex} and Lemma~\ref{brick}, we will be done.
\smallskip

Set
\[
a=\sqrt{\sum_{i<j}(p_ip_j)^2-\binom{n}{2}+1}\,\,\,\, {\rm and}\,\,\,\, x_{ij}=\sqrt{1-(p_ip_j)^2}
\]
for every $i<j$. Let $T_0$ be a regular simplex of side length $a$ with $n$ vertices. Let $S_{ij}$  be a regular simplex of side length $x_{ij}$ with $n-1$ vertices, $1\le i<j\le n$. For the Cartesian product of these simplices,
\[
R=T_0\times\prod_{i<j}S_{ij}\,,
\]
we have
\[
\diam^2(R)=a^2+\sum_{i<j}x_{ij}^2=1\,,
\]
as required.
\smallskip

Let $\pi_0\colon R\longrightarrow T_0$ and $\pi_{ij}\colon R\longrightarrow S_{ij}$
denote the canonical projections. Choose $n$ points, $q_1, \ldots, q_n\in R$
such that
\[
T_0=\{\pi_0(q_1), \ldots, \pi_0(q_n)\}\,,\,\,\,
S_{ij}=\{\pi_{ij}(q_1), \ldots, \pi_{ij}(q_n)\}\,\,\,
\text{ and } \,\,\, \pi_{ij}(q_i)=\pi_{ij}(q_j)\,,
\]
for $1\le i<j\le n$. It remains to check that the
simplex $\{q_1, \ldots, q_n\}$ is congruent to $S$. However, this is obvious, because
\[
(q_kq_\ell)^2=a^2+\sum_{i<j}x_{ij}^2-x_{k\ell}^2=1-x_{k\ell}^2=(p_kp_\ell)^2\,,
\]
for every $1\le k<\ell\le n$.
\hfill $\Box$

\section{\bf Proof of Theorem 5}
Throughout this section, let $d\ge 6$, let $p_0$ denote the origin of $\mathbb{R}^d$, and let  $S=\{p_0,p_1,p_2,p_3\}\subset\mathbb{R}^d$ be the vertex set of a regular tetrahedron of side length $\sqrt2$. Further, let $P\subset\mathbb{R}^d$ denote the $7$-element set consisting of the origin $p_0\in \mathbb{R}^d$ and the (endpoints of the) first $6$ unit coordinate vectors $q_1=(1,0,0,0,0,0,\ldots),$ $q_2=(0,1,0,0,0,0,\ldots),$ $\ldots,$ $q_6=(0,0,0,0,0,1,\ldots).$
Obviously, we have $\diam(S)=\diam(P)=\sqrt2$.
\smallskip

In view of Theorem~\ref{degthm}, in order to establish Theorem 5, it is sufficient to prove that $P$ is $3$-degenerate. That is, we have to show that $\diam(P\cup S)> \sqrt2$. In other words, we have to establish

\begin{clm}\label{claim}
There exist integers $i$ and $j$\, $(1\le i\le 3,\, 1\le j\le 6)$ with $p_iq_j>\sqrt2$.
\end{clm}

The rest of this section is devoted to the proof of this claim.
\smallskip

For $i=1,2,3$, decompose $p_i$ into two components: the orthogonal projection of $p_i$ to the subspace induced by the first 6 coordinate axes and its orthogonal projection to the subspace induced by the remaining coordinate axes. That is, if $p_i=(x_i(1),\ldots,x_i(d))$, let $p_i=p'_i+p''_i$, where
\[
p_i'=(x_i(1),\ldots,x_i(6),0,\ldots,0)\;\;\;{\rm and}\;\;\;
p_i''=(0,\ldots,0,x_i(7),\ldots,x_i(d))\,.
\]
Obviously, we have
\begin{equation}\label{egyenlet1}
|p_i|^2=|p_i'|^2+|p_i''|^2=2\,.
\end{equation}

The proof of Claim~\ref{claim} is indirect. Suppose, for the sake of contradiction, that
\[
\diam\{p_0,p_1,p_2,p_3,q_1,\ldots,q_6\}=\sqrt2\,.
\]
Since $q_j$ and $p_0$ differ only in their $j$th coordinate and $p_iq_j\le p_ip_0$,
the points $p_i$ and $q_j$ lie on the same side of the hyperplane perpendicularly
bisecting the segment $p_0q_j$. That is,
\begin{equation}\label{egyenlet2}
x_i(j)\ge \frac12\;\;\;\text{for every}\;i,j\;\; (1\le i\le 3,\, 1\le j\le 6).
\end{equation}
Hence, we have $|p_i'|^2=\sum_{j=1}^6x_i^2(j)\ge\frac32$ and, by (\ref{egyenlet1}),
\begin{equation}\label{egyenlet3}
|p_i''|^2=|p_i|^2-|p_i'|^2\le\frac12\;\;\;\text{for every}\;i\;\; (1\le i\le 3).
\end{equation}
Moreover, if $i, i'\in\{1,2,3\}$ are distinct, then
\[
\langle p_i, p_{i'}\rangle=\tfrac12\bigl(|p_i|^2+|p_{i'}^2|-|p_i-p_{i'}|^2\bigr)
=\tfrac12(2+2-2)=1\,,
\]
whence~\eqref{egyenlet2} implies
\[
\langle p''_i, p''_{i'}\rangle=1-\sum_{j=1}^6x_i(j)x_{i'}(j)\le-\tfrac12\,.
\]
In view of~\eqref{egyenlet3} it follows that
\[
|p''_1+p''_2+p''_3|^2
=|p''_1|^2+|p''_2|^2+|p''_3|^2
+2\bigl(\langle p''_1, p''_{2}\rangle+\langle p''_1, p''_{3}\rangle
+\langle p''_2, p''_{3}\rangle\bigr)\le -\tfrac32\,,
\]
which is a contradiction. This concludes the proof of Claim~\ref{claim}
and, hence, also the proof of Theorem~5.

\section{Concluding remarks}
\noindent{\bf I. Kneser graphs and hypergraphs.}  Let $d=rn+(k-1)(r-1)$,
where $r,k\ge 2$ are integers. Assign to each $n$-element subset $X\subseteq[d]$
the characteristic vector of $X$. That is, assign to $X$ the point
$p(X)\in\mathbb{R}^d$, whose $i$-th coordinate is
\[
p_i(X)=\left\{
  \begin{array}{ll}
  1 & {\mbox {\rm if}}\;\; i\in X , \\
  0 & {\mbox {\rm if}}\;\; i\not\in X.
  \end{array}
  \right.
\]
Let $P\subseteq\RR^d$ be the set of all points $p(X)$.
We have $|P|=\binom{d}{n}$ and $\diam(P)=\sqrt{2n}$.

\smallskip

For $r=2$, we have $P\subset\mathbb{R}^{2n+k-1}$, and the diameter graph
$H_2(P)$ is called a {\em Kneser graph}. It was conjectured by Kneser~\cite{Kn55}
and proved by Lov\'asz~\cite{Lo78} that $\chi(H_2(P))> k.$ On the other hand, if
$k\le n$, we have $H_3(P)=\emptyset$.

This was generalized to any value of $r$ by Alon, Frankl, and
Lov\'asz~\cite{AlFL86}, who showed that $\chi(H_r(P))>k$, while
$H_{r+1}(P)=\emptyset$, provided that $(k-1)(r-1)<n$. In other words,
the fact that the chromatic number of the $r$-uniform diameter hypergraph of a
point set is high does not imply that the same must hold for its $(r+1)$-uniform
counterpart.

\smallskip

For any integers $r, d\ge 2,$ let $\chi_r(d)$ denote the maximum chromatic number which
an $r$-uniform diameter hypergraph of a point set $P\subseteq\mathbb{R}^d$ can have.

\begin{quest}
Is it true that for every $r\ge 2$, we have $\chi_{r+1}(d)=o(\chi_r(d))$, as $d$ tends to infinity?
\end{quest}

\bigskip

\noindent{\bf II. Relaxations of the diameter-Ramsey property.} Diameter-Ramsey configurations seem to constitute a somewhat peculiar
subclass of the class of all Ramsey configurations. We suggest to classify
all Ramsey configurations $P$ according to the growth rate of the minimum diameter of a point set $R$ with $R\longrightarrow (P)_r$, as $r\rightarrow\infty$.

\begin{dfn}
Given a Ramsey configuration $P$ and an integer $r$, we define
\[
d_P(r)=\inf\,\{\diam(R)\,|\,R\longrightarrow (P)_r\}\,.
\]
\end{dfn}

We have $d_P(r)\ge \diam(P)$, for any Ramsey set $P$ and any integer
$r$, and this holds with equality if and only if for every $\eps>0$ there
exists a configuration $R$ with $R\longrightarrow (P)_r$ and
$\diam(R)\le (1+\eps)\diam(P)$. Certainly, all diameter-Ramsey sets $P$ satisfy
$d_P(r)=\diam(P)$ for all~$r$, but perhaps the configurations with the latter
property form a broader class.

\begin{dfn}
We call a Ramsey set $P$, lying in some Euclidean space,
\begin{enumerate}
\item[(a)] {\it almost diameter-Ramsey} if $d_P(r)=\diam(P)$ holds for all
positive integers $r$;
\item[(b)] {\it diameter-bounded} if there is $C_P>0$ such that
$d_P(r)<C_P$ holds for every positive integer~$r$;
\item[(c)] {\it diameter-unbounded} if $d_P(r)$ tends to infinity, as $r\rightarrow \infty$.
\end{enumerate}
\end{dfn}

We do not know whether there exists any almost diameter-Ramsey configuration
that fails to be diameter-Ramsey. Thus, we would like to ask the following

\begin{quest}
Is it true that every almost diameter-Ramsey set is diameter-Ramsey?
\end{quest}

To establish the diameter-boundedness of certain sets, we may utilize a result of
Matou\v sek and R\"odl~\cite{MaR95}. They showed that, given a
simplex $S$ with circumradius $\varrho$, any number of colors $r$, and
any $\eps>0$, there exists an integer $d$ such that the $d$-dimensional sphere
of radius $\varrho+\eps$ contains a configuration $R$ with
$R\longrightarrow (S)_r$. In particular, this implies the
following

\begin{cor}
Every simplex is diameter-bounded Ramsey.
\end{cor}

Consequently, every diameter-unbounded Ramsey set must be affinely dependent. 
We cannot decide whether there exists any diameter-unbounded Ramsey set, 
but the regular pentagon may serve as a good candidate. K\v{r}\'i\v{z}'s proof 
establishing that the regular pentagon is Ramsey~\cite{Kr91} does not seem to 
imply that it is also diameter-bounded.

\begin{quest}
Is the regular pentagon diameter-unbounded?
\end{quest}

Finally we mention that one can also define these notions for families 
of configurations and ask, e.g., whether they be uniformly diameter-bounded
Ramsey. As an example, we remark that a slight modification of a colouring appearing 
in~\cite{ErGM73} shows that no bounded subset of any Euclidean space 
can simultaneously arrow all triangles whose diameter is $2$ with $8$ colours.
To see this, one may colour each point $x$ with the residue class of $\lfloor 2\|x\|^2\rfloor$
modulo $8$. Given any $K>1$ we set $\xi=\frac{1}{17K^2}$ and consider the isosceles 
triangle with legs of length $1+\xi$ and base of length $2$. Assume for the sake 
of contradiction that there is a monochromatic copy $abc$ of this triangle with 
apex vertex $b$ and with $\|a\|, \|b\|, \|c\|\le K$. Let $m$ denote the mid-point of the 
segment $ac$ and observe that $bm=\sqrt{\xi}$. The triangle inequality yields
\[
	\sqrt{\xi}=\|b-m\|\ge \big|\|b\|-\|m\|\big|
\]
and, hence, we have 
\[
	\sqrt{\xi}\cdot \bigl(\|b\|+\|m\|\bigr)\ge \big|\|b\|^2-\|m\|^2\big|\,.
\]
Multiplying by $4$, and applying triangle inequality to the left-hand side and the 
parallelogram law to the right-hand side we infer
\begin{align*}
	2\sqrt{\xi}\cdot \bigl(\|a\|+2\|b\|+\|c\|\bigr) & \ge \big| 4\|b\|^2-\|a+c\|^2\big| \\
	& = \big| 4\|b\|^2-2\|a\|^2-2\|c\|^2+\|a-c\|^2\big| \\
	& = \big|4+\bigl(2\|b\|^2-2\|a\|^2\bigr) +\bigl(2\|b\|^2-2\|c\|^2\bigr)\big|\,,
\end{align*}
which due to 
$\lfloor 2\|a\|^2\rfloor\equiv \lfloor 2\|b\|^2\rfloor\equiv\lfloor 2\|c\|^2\rfloor\pmod{8}$
leads to $8K\sqrt{\xi}\ge 2$, contrary to our choice of $\xi$.

\begin{rem}
While revising this article, we learned from Nora Frankl about some 
progress regarding Question~\ref{qu:obtuse} obtained jointly with Jan Corsten~\cite{CF17}. 
They proved that the bound of $150^\circ$ appearing in Theorem~4 above can be lowered 
to $135^\circ$. Their elegant proof involves the spherical colouring and Jung's inequality. 
\end{rem}

\end{document}